\documentclass[reqno]{amsart}
\usepackage{amsmath,amssymb,colordvi}
\usepackage{multicol}
\usepackage{color,bbm}

 \newtheorem{theorem}{Theorem}[section]
 
 \newtheorem{lemma}[theorem]{Lemma}

 \newtheorem{definition}[theorem]{Definition}
 \newtheorem{remark}[theorem]{Remark}
  
 \newtheorem{assumption}[theorem]{Assumption}
 \numberwithin{equation}{section}

 \newcommand{\RR}{\mathbb{R}}
 \newcommand{\NN}{\mathbb{N}}
 \newcommand{\CC}{\mathbb{C}}

 \newcommand{\abs}[1]{\left\vert#1\right\vert}

 \newcommand{\norm}[1]{\left\Vert#1\right\Vert}

  \newcommand{\Om}{\Omega}
 \def\pOm{\partial\Omega}

 \usepackage[dvips,bottom=1.4in,right=1in,top=1in, left=1in]{geometry}

\def\D{{\mathbb{D}}}
\def\RR{{\mathbb{R}}}
\def\NN{{\mathbb{N}}}

\def\CC{{\mathbb{C}}}

\def\Om{\Omega}

\def\pOm{\partial\Omega}

\begin{document}

\title[Memory approximate controllability]{Memory approximate controllability properties for higher order Hilfer time fractional evolution equations }

\author{Ernes Aragones}
\address{%
E. Aragones, University of Puerto Rico, Carolina Campus\\
Department of Natural Sciences \\
 PO Box 4800, Carolina PR 00984-4800 (USA)}
\email{ernest.aragones@upr.edu}

\author{Valentin Keyantuo}
\address{%
V. Keyantuo, University of Puerto Rico, Rio Piedras Campus\\
Faculty of Natural Sciences, Department of Mathematics \\
 17 University AVE. STE 1701, San Juan PR 00925-2537 (USA)}
 \email{valentin.keyantuo1@upr.edu}

\author{Mahamadi Warma}
\address{M. Warma,  Department of Mathematical Sciences and the Center for Mathematics and Artificial Intelligence, George Mason University.  Fairfax VA 22030 (USA)}
\email{mwarma@gmu.edu}

\thanks{The work of the authors is partially supported by the US Army Research Office (ARO) under Award NO: W911NF-20-1-0115}



\begin{abstract}
In this paper we study the approximate controllability of fractional partial differential equations associated with the so-called Hilfer type time fractional derivative and a non-negative selfadjoint operator $A$ with a compact resolvent on $L^2(\Omega)$, where $\Omega\subset\RR^N$ ($N\geq 1$) is an open set. More precisely, we show that if $0\le\nu\le 1$, $1<\mu\le 2$ and $\Omega\subset\RR^N$ is an open set, then the system 
\begin{equation*}
\begin{cases}
\D^{\mu,\nu}_tu+Au=f\chi_{\omega}\;\;&\mbox{ in }\;\Omega\times(0,T),\\
(I_t^{(1-\nu)(2-\mu)}u)(\cdot,0)=u_0 &\mbox{ in }\;\Omega,\\
(\partial_tI_t^{(1-\nu)(2-\mu)}u)(\cdot,0)=u_1 &\mbox{ in }\;\Omega,
\end{cases}
\end{equation*}
is memory approximately controllable for any $T>0$, $u_0\in D(A^{1/\mu})$, $u_1\in L^2(\Omega)$ and any non-empty open set $\omega\subset\Omega$. The same result holds for every $u_0\in D(A^{1/2})$ and $u_1\in L^2(\Omega)$.
\end{abstract}

\maketitle

\section{Introduction}
Fractional order differential equations have attracted interest from researchers due to their suitability as modeling tools for several phenomena in science and technology. Specifically, phenomena with memory in areas such as viscoelasticity, rheology, anomalous diffusion and many more have been found to be better modeled through the use of fractional differential equations. It is also true that in the space variables, the use of fractional order operators is witnessing  an almost explosive development.

Fractional order equations are not new, and the idea of fractional derivative goes back to the birth of differential calculus with a question of de L'H\^opital to Leibniz (\cite{Po99}, \cite{Mi-Ro}, \cite{Ma97}, \cite{Ma}). Riemann and Liouville worked on fractional integrals and derivatives. We should also mention Abel's equation in the modeling of the tautochrone (\cite{podlubny2017niels}, \cite{gorenflo1991abel}). Subsequently, H. Weyl and several other mathematicians contributed to the development of fractional calculus. The theory of function spaces clearly makes contact with  fractional calculus.  Indeed, scales of function spaces such as fractional order Sobolev spaces, Besov spaces and Triebel-Lizorkin spaces are used to achieve more precise regularity results in partial differential equations and variational calculus. Fractional calculus is mathematically intimately related to the theory of integral equations.

Concerning evolution equations with fractional order derivatives, the Riemann-Liouville derivative is the oldest. A modification thereof was proposed in the last century by F. Caputo (\cite{Cap}). This modified version, nowadays known as the Caputo fractional derivative, has the advantage that in its use, initial conditions are similar to what is known for ordinary derivatives. Yet, the Riemann-Liouville fractional derivative which necessitates that the initial conditions be given in integral form, has been found to be more suitable for some applications (\cite{Pod}).

Most recently R. Hilfer (\cite{Hil-2000}) has proposed a version of fractional derivative which interpolates between the Riemann-Liouville and the Caputo derivatives. In the present paper, we shall investigate approximate controllability of evolution problems involving this Hilfer fractional derivative. For a locally integrable function $f:\, [0,\,\infty)\longrightarrow X$ where $X$ is a Banach space, and $\alpha>0,$ the Riemann-Liouville fractional integral of order $\alpha$ is defined by
\begin{equation}\label{IT}
(I_t^\alpha f)(t)=\int_0^t \frac{(t-s)^{\alpha-1}}{\Gamma(\alpha)} f(s)\;ds,\,\, t>0
\end{equation}
where $\Gamma$ denotes the usual Euler-Gamma function. Given $T>0$,  the right-hand Riemann-Liouville fractional integral of order $\alpha$ on $(0,T)$ is defined by
\begin{equation}\label{ItT}
(I_{t,T}^\alpha f)(t):= \frac{1}{\Gamma(\alpha)}\int_{t}^{T}(s-t)^{\alpha-1}f(s)\;ds,\quad (t<T).
\end{equation}
The right-hand Riemann-Liouville fractional derivative of order $\alpha$ with $n-1\le \alpha<n$, $n\in\mathbb N$ is given by
\begin{equation}\label{LRLFD}
D_{t,T}^\alpha f(t)=(-1)^n\frac{d^n}{dt^n}I_{t,T}^{n-\alpha}f(t),\quad (t<T).
\end{equation}

Then,  for $f$ satisfying appropriate conditions, the Hilfer fractional derivative of order $(\mu, \nu)$ is given by
\begin{equation}\label{Hilfer0001}
\mathbb D_t^{\mu,\nu}f(t)=I_t^{\nu(1-\mu)}\frac{d}{dt}I_t^{(1-\nu)(1-\mu)}f(t),\, \, t>0.
\end{equation}
In \eqref{Hilfer0001},  it is assumed that $0<\mu\le 1$ and $0\le \nu\le 1.$ This case is studied in \cite{AKW}.

We shall be concerned here with the higher order case, that is, $1<\mu\le 2$ and $0\le \nu\le 1.$ Under these conditions, we set 
\begin{equation}\label{Hilfer0002}
\mathbb D_t^{\mu,\nu}f(t)=I_t^{\nu(2-\mu)}\frac{d^2}{dt^2}I_t^{(1-\nu)(2-\mu)}f(t),\, \, t>0.
\end{equation}
It is worthwhile to observe that when $\mu=2$,  \eqref{Hilfer0002} reduces to the second order derivative while for $\nu=0,\, 1<\mu\le 2$, it reduces to the Riemann-Liouville fractional derivative of order $\mu$. Moreover, for $\nu=1$, we find the Caputo fractional derivative of order $\mu$.
We should also point out that in the limiting case where $\mu=1$ and $\nu=0$, we obtain $\mathbb D_t^{1,0}f(t)=f^\prime(t).$ So,  naturally, the Hilfer derivative proposed in \eqref{Hilfer0002} interpolates between the first and second derivatives.

Let $\Omega\subset\RR^N$ ($N\ge 1$) be an open set with boundary $\pOm$.   The regularity needed on $\Omega$ will depend on the operator $A$ under consideration so that our assumptions (see Assumption \ref{assump}) are satisfied.
The main concern of the present paper is to study the controllability properties of a class of fractional (possible space-time) differential equations involving the Hilfer time-fractional derivative. More precisely, we consider the following initial  value problem:
\begin{equation}\label{main-EQ}
\begin{cases}
\D^{\mu,\nu}_tu+Au=f\chi_{\omega}\;\;&\mbox{ in }\;\Omega\times(0,T),\\
(I_t^{(1-\nu)(2-\mu)}u)(\cdot,0)=u_0 &\mbox{ in }\;\Omega,\\
(\partial_tI_t^{(1-\nu)(2-\mu)}u)(\cdot,0)=u_1 &\mbox{ in }\;\Omega,
\end{cases}
\end{equation}
where $T>0$ and $1<\mu\le 2,\, 0\le\nu\le 1$ are real numbers, $\mathbb D_t^{\mu,\nu} u$ denotes the Hilfer time fractional derivative of order $({\mu,\nu})$ of the function $u$ given in \eqref{Hilfer0002} and the operator $A$ is a non-negative selfadjoint operator on $L^2(\Omega)$ with compact resolvent.

 We are interested in the approximate controllability of \eqref{main-EQ}. It has been shown that  for fractional in-time equations, exact or null controllability is in general not achievable (\cite{EZ}). For this reason, approximate controllability is a viable substitute. Approximate controllability for fractional evolution equations has been studied recently by many authors (\cite{FY}, \cite{Ke-wa}, \cite{LiuLi2015}, \cite{War-SIAM}), both for the Caputo fractional derivative of order $\mu$ with $0<\mu<1$ (\cite{FY,War-SIAM}) and $1<\mu<2$ (see \cite{Ke-wa}), and the Riemann-Liouville fractional derivative 
 for $0<\mu<1$ (\cite{LiuLi2015}) . 
 
 In the recent paper \cite{AKW}, the authors introduced the notion of mean approximate controllability and established its validity for the analogue of equation \eqref{main-EQ} when $0<\mu<1$, for a large class of operators $A$.  The notion of memory controllable introduced in the present paper is the same notion as the mean conrollability studied in \cite{AKW}.  We think that memory controllability is more suitable.
We will consider Problem \eqref{main-EQ} in the space $L^2(\Omega)$ where $\Omega$ is an open subset of $\mathbb R^N$ with boundary $\partial\Omega$, and $A$ is a nonnegative self-adjoint operator with compact resolvent. More concretely, the operator $A$ can be an elliptic differential operator with Dirichlet boundary conditions. Another interesting example is the fractional Laplace operator with null exterior Dirichlet condition or null exterior nonlocal Robin type conditions (see e.g.  \cite{BWZ,LR-WA,War-N} for the definition of these operators). 
 
 The rest of the paper is organized as follows. In Section 2, we give some preliminary material on  fractional derivatives, the Mittag-Leffler functions and the associated Laplace transforms.  We state the main results (the unique continuation principle of solutions to the associated adjoint equation and the memory approximate controllability of \eqref{main-EQ}). 
 We also introduce the solution families for \eqref{main-EQ} under the assumptions imposed on the operator $A$. These are in fact resolvent families. We do not include here the general theory of resolvent families and merely give the corresponding representation under our hypotheses.  The special  case of the Caputo fractional derivative of order $1<\mu<2$ appears in the references \cite{AGKW,KLW-2}. Still, for the Caputo fractional derivative in the case where $A$ is a positive self-adjoint operator with compact resolvent, the representation is given in \cite{Ke-wa} and \cite{KLW-2}. Then follows the proof of well-posedness of \eqref{main-EQ}, along with some estimates on the solutions. The adjoint system whose well-posedness is needed in the proof of our main results is studied in this section. Sections \ref{sec3} and \ref{sec4} are devoted to the proof of the main results.   
 
 

\section{Preliminaries, main results, existence and regularity of solutions}

In this section we give some well-known results that are used throughout the paper,  state our main controllability results, and
prove some existence and regularity results of the  systems under study.

\subsection{Preliminaries}
Let $X$ be a Banach space. We have the following useful result. It can be proved by using the Laplace transform. Here, we give a direct proof.

\begin{lemma}
Let $f,g:[0,\infty)\to X$ be locally integrable and $I_t^\alpha$ be the operator defined in \eqref{IT}.
Then,  for every $t>0$ and $\alpha>0$,
\begin{align}\label{Con-Int}
((I^\alpha_t f)\ast g)(t)=(f\ast(I^\alpha_t g))(t).
\end{align}
\end{lemma}

\begin{proof}
Using the definition and a change of variables,  a simple calculation gives 
\begin{align*}
((I^\alpha_t f)\ast g)(t)&=\int_0^tg(t-s)(I^\alpha_sf)(s)\,ds\\
&=\dfrac{1}{\Gamma(\alpha)}\int_0^t\int_0^s\,g(t-s)(s-\tau)^{\alpha-1}f(\tau)\,d\tau ds\\&=\dfrac{1}{\Gamma(\alpha)}\int_0^t\int_\tau^t\,g(t-s)(s-\tau)^{\alpha-1}f(\tau)\,ds d\tau  \;\;\;\;(\mbox{changes of variables})\\&=\dfrac{1}{\Gamma(\alpha)}\int_0^t\int_0^{t-\tau}\,g(\eta)(t-\tau-\eta)^{\alpha-1}f(\tau)\,d\eta\, d\tau \;\;\;\;(\eta=t-s)\\&=\int_0^t\,(I^\alpha_{t-\tau}g)(t-\tau)f(\tau)\,d\tau\\&=((I^\alpha_t g)\ast f)(t)
\end{align*}
and the proof is finished.
\end{proof}

We have the following property for power functions.

\begin{lemma}\label{int-frac}
Let  $\alpha>0$ and $\beta>-1$. Then,
	\begin{align*}
	 I_t^\alpha(t^\beta)=\dfrac{\Gamma(\beta+1)}{\Gamma(\alpha+\beta+1)}t^{\alpha+\beta},\qquad t>0.
	\end{align*}		
\end{lemma}

The integral operators $I_t^\alpha$ and $I^\alpha_{t,T}$ satisfy the following integration by parts formula. The proof is contained  in \cite[Corollary of Theorem 3.5]{Ma} (see also e.g. \cite{Agr,AlTo, THS10} and the references therein).

\begin{lemma} 
Let  $\alpha>0$, $\varphi\in L^p(0,T)$ and $\psi\in L^q(0,T)$ with $p\geq 1,\;q\geq 1$ and $\frac 1p+\frac 1q\leq 1+\alpha$, but $p\neq 1$, $q\neq 1$ in the case $\frac 1p+\frac 1q= 1+\alpha$. Then,
\begin{align}\label{IPF}
\int_0^T\varphi(t)(I^\alpha_t\psi)(t)dt=\int_0^T\psi(t)(I^\alpha_{t,T}\varphi)(t)\,dt.
\end{align}
\end{lemma}


Also of interest is the right-hand Hilfer time-fractional derivative of order $(\mu,\nu)$ ($0\le  \nu\le 1$, $1<\mu \le 2$) given by
\begin{align}\label{D-R}
\mathbb{D}_{t,T}^{\mu,\nu} u(t):=& - I_{t,T}^{\nu(2-\mu)} \frac{d^2}{dt^2}\left( I_{t,T}^{(1-\nu)(2-\mu)} u\right)(t).
\end{align}
The right-hand derivatives and integrals are introduced above  since they are needed for the integration by parts formula. Indeed, we have the following  integration by parts formula  (see e.g. \cite{Agr,AlTo,THS10}).
 
\begin{lemma} 
Let $0\le  \nu\le 1$ and $1<\mu \le 2$.  
Then, 
\begin{align}\label{IP01} 
\int_0^Tv(t)  \mathbb D_{t}^{\mu,\nu} u(t)\;dt=&\int_0^T u(t)\mathbb D_{t,T}^{\mu,1-\nu} v(t)\;dt\notag\\
&+\left[ I_{t,T}^{(1-\nu)(2-\mu)}u(t)D_{t,T}^{1-\nu(2-\mu)}v(t)+\dfrac{d}{dt}I_{t,T}^{(1-\nu)(2-\mu)}u(t)I_{t}^{\nu(2-\mu)}v(t)\right]_{t=0}^{t=T},
\end{align}
 provided that the left and right-hand side expressions make sense.
 \end{lemma} 
 
 \begin{proof}
  Since 
 $$\mathbb D_t^{\mu,\nu}u(t)=I_t^{\nu(2-\mu)}\dfrac{d^2}{dt^2}\left( I_t^{(1-\nu)(2-\mu)}u(t)\right)  \mbox{ with } 0<\nu(2-\mu),(1-\nu)(2-\mu)<1,$$
 using Equation \eqref{IPF} we have that,
  \begin{align}\label{IP02}
\int_0^Tv(t)  \mathbb D_{t}^{\mu,\nu} u(t)\;dt=\int_0^T \dfrac{d^2}{dt^2}\left( I_t^{(1-\nu)(2-\mu)}u\right)(t)(I_{t,T}^{\nu(2-\mu)}v)(t)\,dt.
\end{align}  
   Now, applying classical integration by parts twice in the right hand side of \eqref{IP02} gives
  \begin{align*}
\int_0^Tv(t)  \mathbb D_{t}^{\mu,\nu} u(t)\;dt=&\left[\dfrac{d}{dt}(I_t^{(1-\nu)(2-\mu)}u)(t)(I_{t,T}^{\nu(2-\mu)}v)(t)-(I_t^{(1-\nu)(2-\mu)}u)(t)\dfrac{d}{dt}(I_{t,T}^{\nu(2-\mu)}v)(t) \right]_{t=0}^{t=T}\\ 
&+\int_0^T\dfrac{d^2}{dt^2}(I_t^{\nu(2-\mu)}v)(t)(I_t^{(1-\nu)(2-\mu)}u)(t)\;dt.
\end{align*}  
Using again \eqref{IPF} in the last integral, we get  \eqref{IP01}  where we recall that 
  $$I_{t,T}^{(1-\nu)(2-\mu)}\dfrac{d^2}{dt^2}(I_t^{\nu(2-\mu)}v)(t)=\mathbb D_{t,T}^{\mu,1-\nu}v(t)\mbox{ and } -\dfrac{d}{dt}(I^{\nu(2-\mu)}_{t,T}v)(t)=D_{t,T}^{1-\nu(2-\mu)}v(t).$$
  The proof is finished.
 \end{proof}

 Special cases related to the Caputo and Riemann-Liouville fractional derivatives are easily obtained from the above formula \eqref{IP01}.
 
The Mittag-Leffler function with two parameters is defined as follows:
\begin{align}\label{mm}
E_{\alpha, \beta}(z) := \sum_{n=0}^{\infty}\frac{z^n}{\Gamma(\alpha n + \beta)},\; \;\alpha>0,\;\beta \in\CC, \quad z \in \CC.
\end{align}
It is well-known that $E_{\alpha,\beta}(z)$ is an entire function. This is so even if we allow the parameter set to include $\alpha\in\mathbb C$ with $\rm{Re}(\alpha)>0.$

The following estimate of the Mittag-Leffler function will be useful. Let $0<\alpha\le2$, $\beta\in\RR$ and $\kappa$ be such that $\frac{\alpha\pi}{2}<\kappa<\min\{\pi,\alpha\pi\}$. Then,  there is a constant $C=C(\alpha,\beta,\kappa)>0$ such that
\begin{equation}\label{Est-MLF}
|E_{\alpha,\beta}(z)|\le \frac{C}{1+|z|},\;\;\;\kappa\le |\mbox{arg}(z)|\le \pi.
\end{equation}
The Mittage-Leffler function with two parameters given in \eqref{mm} satisfies the following recurrence relation (see e.g. \cite[Formula (4.2.4), p.57]{GKMR2020} or \cite[Theorem 5.1]{HMS}):
\begin{align}\label{MLD}
E_{\alpha,\beta}=\beta E_{\alpha,\beta+1}(z)+\alpha z\dfrac{d}{dz}E_{\alpha,\beta+1}(z).
\end{align}

We have the following  useful estimates of the Mittag-Leffler functions that can be verified by using \eqref{Est-MLF}. We refer to \cite{Ke-wa} for the full proof.

\begin{lemma}\label{lem-INE} Let $1<\alpha<2$ and $\beta>0$. Then the following assertion hold.
\begin{enumerate}
\item Let $0\le \nu \le 1$, $0<\gamma<\alpha$ and $\lambda>0$. Then there is a constant $C>0$ such that for every $t>0$,
\begin{align}\label{Est-ML2}
\abs{\lambda^\nu t^\gamma E_{\alpha,\beta}(-\lambda t^\alpha)}\le Ct^{\gamma-\alpha\nu}.
\end{align}
\item Let $0\le\gamma\le 1$ and $\lambda>0$. Then there is a constant $C>0$ such that for every $t>0$,
\begin{align}\label{Est-ML3}
\abs{\lambda^{1-\gamma}t^{\alpha-2}E_{\alpha,\beta}(-\lambda t^\alpha)}\le Ct^{\alpha\gamma-2}.
\end{align}
\end{enumerate}
\end{lemma}

The Laplace transform  of the Mittag-Leffler function is given by the following important relation:
\begin{equation}\label{lap-ml}
\int_0^{\infty} e^{-\lambda t} t^{\alpha k + \beta-1}
E_{\alpha,\beta}^{(k)}(\pm \gamma t^{\alpha})dt = \frac{k!
\lambda^{\alpha-\beta}}{(\lambda^{\alpha} \mp \gamma)^{k+1}},
\quad \mbox{Re}(\lambda)> |\gamma|^{1/\alpha}.
\end{equation}
Here, $k\in\mathbb{N}\cup\{0\}$, $\alpha>0$ and $\gamma\in\mathbb{R}.$ 

The following formulas can be proved by using the Laplace transform (see e.g. \cite{Ke-wa}).

\begin{lemma}\label{Der-MLF}
Let $\alpha>0$, $\beta>0$, $\lambda>0$, $t>0$ and $m\in\NN$. Then,
\begin{align}\label{Der-ML-1}
\dfrac{d^m}{dt^m}\left[ E_{\alpha,1}(-\lambda t^\alpha)\right] =-\lambda t^{\alpha-m}E_{\alpha,\alpha-m+1}(-\lambda t^\alpha), 
\end{align}
\begin{align}\label{Der-ML-2}
\dfrac{d}{dt}\left[ tE_{\alpha,2}(-\lambda t^\alpha)\right] =E_{\alpha,1}(-\lambda t^\alpha),
\end{align}
\begin{align}\label{Der-ML-3}
\dfrac{d}{dt}\left[ t^{\alpha-1}E_{\alpha,\alpha}(-\lambda t^\alpha)\right] =t^{\alpha-2}E_{\alpha,\alpha-1}(-\lambda t^\alpha),
\end{align}
\begin{align}\label{EE-11}
\int_0^tE_{\alpha,1}(-\lambda \tau^\alpha)\,d\tau=tE_{\alpha,2}(-\lambda t^\alpha).
\end{align}
\end{lemma}
%

 For more details on fractional derivatives, integrals and the  Mittag-Leffler functions  we refer to \cite{Agr, Ba01,Go-Ma97,Ma97,Go-Ma00,Mi-Ro,Po99} and the references therein.

\subsection{Main results} \label{main-res}

In this section we state the main controllability results of the paper.  The notion of solutions to the  systems under study will be introduced in the next section.  First,  we make the following assumption on the operator $A$.

\begin{assumption}\label{assump}
We assume that $A:D(A)\subset L^2(\Omega)\to L^2(\Omega)$ is a nonnegative selfadjoint operator which is invertible and has  a compact resolvent.  
\end{assumption}

We refer to \cite[Section 2.2]{AKW} for several examples of operators that enter in our framework.

\begin{remark}
{\em
As a consequence of Assumption \ref{assump} we have the following.
\begin{enumerate}
\item By the spectral theorem, one can define the powers  $A^\beta$ of the operator $A$ for any $\beta\in\mathbb R$.
\item The operator $A$ is given by a bilinear,  symmetric, coercive and closed form $\mathcal E_A:V_{1/2}\times V_{1/2}\to \mathbb R$ given by 
$$\mathcal E_A(u,v)=(A^{1/2}u,A^{1/2}v)_{L^2(\Om)} \mbox{ for all } u,v\in V_{1/2}:=D(A^{1/2}).$$

\item For every $\gamma\ge 0$, we denote  $V_\gamma:=D(A^\gamma)$ and $V_{-\gamma}$ the dual of $V_{\gamma}$ with respect to the pivot space $L^2(\Omega)$ so that we have the  continuous and dense embeddings $V_\gamma \hookrightarrow L^2(\Omega)\hookrightarrow V_{-\gamma}$.   Notice that $V_{-\gamma}=D(A^{-\gamma})$.

\item We denote by $(\varphi_n)$ the orthonormal basis of eigenfunctionss of $A$ associated with the eigenvalues $(\lambda_n)$. Then,  $0<\lambda_1\le\lambda_2\le\cdots\le \lambda_n\le\cdots$ and $\lim_{n\to\infty}\lambda_n=+\infty$.
\end{enumerate}
}
\end{remark}

Next, we introduce the notion of memory approximate controllability.


\begin{definition} \label{def-211} 
Let $0\le \nu\le 1$, $1<\mu\le 2$, $\gamma:=1/\mu$ and $p>1/(\mu-1)$ (or $\gamma=1/2$ and $p>2/\mu$). 
We say that the system \eqref{main-EQ} is memory approximately controllable in time $T>0$, 
if for every $(u_0,u_1)$,  $(v,w)\in V_\gamma\times L^2(\Omega)$ and $\varepsilon>0$, there exists a control function $f\in L^p((0,T);L^2(\omega))$ such that the corresponding weak solution $u$ of \eqref{main-EQ} (see Definition \ref{def-214}) satisfies
\begin{align*}
\norm{\left(I_t^{(1-\nu)(2-\mu)}u\right)(\cdot,T)-v}_{V_\gamma}+\norm{\left(\partial_t I_t^{(1-\nu)(2-\mu)}u\right)(\cdot,T)-w}_{L^2(\Omega)}\le\varepsilon.
\end{align*}
This is equivalent to requiring that the set
$$\left\{\left[\left(I_t^{(1-\nu)(2-\mu)}u\right)(\cdot,T), \left(\partial_t I_t^{(1-\nu)(2-\mu)}u\right)(\cdot,T)\right]:\; f\in L^p((0,T);L^2(\omega))\right\}$$
is dense in $V_\gamma\times L^2(\Omega)$.
\end{definition}

The following theorem is our first main result.

\begin{theorem}\label{second-Theo} 
Let $0\le \nu\le 1$, $1<\mu\le 2$ and $\omega\subset\Omega$ an arbitrary non-empty open set. Assume that the operator $A$ has the unique continuation property in the sense that,
\begin{align}\label{ucp}
\mbox{ if }\lambda>0,\;\; \varphi\in D(A),\;\; A\varphi=\lambda\varphi\;\mbox{ and } \varphi=0\mbox{ in }\omega, \mbox{ then } \varphi=0 \mbox{ in }\Omega.
 \end{align}
Then, the system \eqref{main-EQ} is memory approximately controllable in any time $T>0$.
\end{theorem}

Using the integration by parts formula \eqref{IP01} it is straightforward to show that the following system
\begin{equation}\label{ACP-Dual}
\begin{cases}
\mathbb D_{t,T}^{\mu,1-\nu} v +Av=0\;\;&\mbox{ in }\; \Omega\times (0,T),\\
\left(I_{t,T}^{\nu(2-\mu)}v\right)(\cdot,T)=v_0   \;&\mbox{ in }\;\Omega, \\
\left( D_{t,T}^{1-\nu(2-\mu)}v\right)(\cdot,T)=v_1    \;&\mbox{ in }\;\Omega,
\end{cases}
\end{equation}
can be viewed as the dual system associated with \eqref{main-EQ}.

Next, we show that under the assumption that $A$ has the unique continuation property, the adjoint system \eqref{ACP-Dual} satisfies the  unique continuation principle which is our second main result.

\begin{theorem}\label{pro-uni-con}
Let $0\le  \nu\le 1$, $1<\mu \le 2$, $\gamma:=1/\mu$ or $\gamma=1/2$, $v_0\in V_\gamma$, $v_1\in L^2(\Omega)$ and let $\omega\subset\Omega$ be an arbitrary non-empty open set. Assume that $A$ has the unique continuation property in the sense of \eqref{ucp}.
Let $v$ be the unique weak solution of \eqref{ACP-Dual}. If $v=0$ in $\omega\times (0,T)$, then $v=0$ in $\Omega\times (0,T)$.
\end{theorem}

We conclude this section with the following observation.

\begin{remark}
{\em As we shall see in the proof of the main results, it turns out that memory approximately controllable is equivalent to the unique continuation principle of solutions to the dual system \eqref{ACP-Dual}, that is, to Theorem \ref{pro-uni-con}.
}
\end{remark}

\subsection{Well-posedness and representation of solutions to Equation \eqref{main-EQ}}

Our notion of weak solutions to the system \eqref{main-EQ} is as follows.

\begin{definition} \label{def-214}
Let  $0\le  \nu\le 1$ and $1<\mu \le 2$.
Let $1/2\le \gamma\le 1$, $p\ge 1$,  $u_0\in V_\gamma$, $u_1\in L^2(\Omega)$ and $f\in L^p((0,T);L^2(\Omega))$. A function $u$ is said to be a weak solution of \eqref{main-EQ}, if for every $T>0$, the following properties hold:
\begin{itemize}
\item Regularity:
\begin{equation}\label{regu}
\begin{cases}
u\in C((0,T];V_\gamma), \\
 I_t^{(1-\nu)(2-\mu)}u\in C([0,T];V_\gamma)\cap C^1([0,T];L^2(\Omega)),\\
\mathbb D_t^{\mu,\nu} u\in C((0,T];V_{-\gamma}),
\end{cases}
\end{equation}
\item Initial Conditions: 
\begin{align}\label{Init-Cond}
 \left(I_t^{(1-\nu)(2-\mu)}u\right)(\cdot,0)=u_0\;\mbox{ and }\; \partial_t I_t^{(1-\nu)(2-\mu)}u(\cdot,0)=u_1, \; \mbox{a.e.\,in}\;\Omega.
\end{align}
\item Variational  identity: for every $\varphi\in  V$ and a.e. $t\in (0,T)$, 
\begin{align}\label{Var-I}
\langle \mathbb D_t^{\mu,\nu} u(\cdot,t),\varphi\rangle_{V_{-\gamma},V_\gamma}+\mathcal E_A(u(\cdot,t),\varphi)=(f(\cdot,t),\varphi)_{L^2(\Om)}.
\end{align}
\end{itemize} 
\end{definition}

\begin{remark}
{\em 
Observe that if $1<\mu\le 2$ and $\gamma=1/\mu$ (hence $1/2\le \gamma<1$) we have  the continuous embedding $V_\gamma\hookrightarrow V_{1/2}$ so that $\mathcal E_A(u,v)$ is well defined for every $u,v\in V_\gamma$.
}
\end{remark}


Recall that $(\lambda_n)$, $0<\lambda_1\le \lambda_2\le \cdots\le \lambda_n\le\cdots$ are the eigenvalues of $A$ and $(\varphi_n)$ the orthonormal basis of eigenfunctions associated with the eigenvalues $(\lambda_n)$. We have the following result on existence, uniqueness and representation of solutions.

\begin{theorem}\label{theo-weak} 
Let $0\le  \nu\le 1$, $1<\mu \le 2$,  $\gamma:=1/\mu$ and $p>1/(\mu-1)$ (or $\gamma=1/2$ and $p>2/\mu$),  and $\beta:=(2-\mu)(1-\nu)$. Then,  for every $u_0\in V_\gamma$, $u_1\in L^2(\Omega)$ and $f\in L^p((0,T);L^2(\Omega))$, the system \eqref{main-EQ} has a unique weak solution $u$ given by
\begin{align}\label{Weak-Sol}
u(\cdot,t)=&\sum_{n=1}^\infty(u_0,\varphi_n)_{L^2(\Om)}t^{-\beta}E_{\mu,1-\beta}(-\lambda_nt^\mu)\varphi_n+\sum_{n=1}^\infty(u_1,\varphi_n)_{L^2(\Om)}t^{1-\beta}E_{\mu,2-\beta}(-\lambda_nt^\mu)\varphi_n \notag\\&+\sum_{n=1}^\infty\left( \int_0^t\,(f(\cdot,\tau),\varphi_n)_{L^2(\Om)}(t-\tau)^{\mu-1}E_{\mu,\mu}(-\lambda_n(t-\tau)^\mu)d\tau \right) \varphi_n.
\end{align}

\end{theorem}

\begin{proof} 
Let $u_{0,n}:=(u_0,\varphi_n)_{L^2(\Omega)}$, $u_{1,n}:=(u_1,\varphi_n)_{L^2(\Omega)}$,  $u_{n}(t):=(u(t),\varphi_n)_{L^2(\Omega)}$ and $f_n(t):=(f(\cdot,t),\varphi_n)_{L^2(\Omega)}$. 
We proceed in several steps.  We consider the case $\gamma=1/\mu$ and $p>1/(\mu-1)$. The case $\gamma=1/2$ and $p>2/\mu$ follows similarly with the appropriate modifications. 

{\bf Step 1}: Multiplying the first equation in \eqref{main-EQ} by $\varphi_n$ and integrating over $\Omega$, we get that $u_n(t)$ is the solution of the following fractional ordinary differential equation:
\begin{equation}\label{ODE}
\begin{cases}
\D^{\mu,\nu}_tu_n(t)+\lambda_nu_n(t)=f_n(t)\;\;\mbox{ in }\;(0,T),\\
(I_t^{(1-\nu)(2-\mu)}u_n)(0)=u_{0,n},\\
(\partial_tI_t^{(1-\nu)(2-\mu)}u)(0)=u_{1,n}.
\end{cases}
\end{equation}
Using the method of the Laplace transform we have that the unique solution of \eqref{ODE} is given for every $n\in\mathbb N$ by
\begin{align*}
u_n(t)=t^{-\beta}E_{\mu,1-\beta}(-\lambda_nt^\mu)u_{0,n} +t^{1-\beta}E_{\mu,2-\beta}(-\lambda_nt^\mu)u_{1,n} +\int_0^t(t-\tau)^{\mu-1}E_{\mu,\mu}(-\lambda_n(t-\tau)^\mu)f_n(\tau)\;d\tau.
\end{align*}
This suggests that a formal solution of \eqref{main-EQ} is given by \eqref{Weak-Sol}. We show that indeed,  \eqref{Weak-Sol} satisfies all the assumptions in Definition \ref{def-214}.
In order to simplify the notations,  for $v\in L^2(\Om)$ we define the operator families
\begin{equation*}
S_1(t)v:=\sum_{n=1}^\infty(v,\varphi_n)_{L^2(\Om)}t^{-(1-\nu)(2-\mu)}E_{\mu,\mu-1+\nu(2-\mu)}(-\lambda_nt^\mu)\varphi_n,\, \, t>0,
\end{equation*}
\begin{equation*}
S_2(t)v:=\sum_{n=1}^\infty(v,\varphi_n)_{L^2(\Om)}t^{\mu-1+\nu(2-\mu)}E_{\mu,\mu+\nu(2-\mu)}(-\lambda_nt^\mu)\varphi_n, \, \, t>0,
\end{equation*}
and
\begin{equation*}
S_3(t)v=\sum_{n=1}^\infty (v, \varphi_n)_{L^2(\Om)} t^{\mu-1}E_{\mu,\mu}(-\lambda t^\mu)\varphi_n,\, \, t>0.
\end{equation*}
Then,  \eqref{Weak-Sol} can be rewritten as
\begin{align*}
u(\cdot,t)=S_1(t)u_0+S_2(t)u_1+(S_3\ast f)(t).
\end{align*}

{\bf Step 2}:  We show that $u\in C((0,T];V_{\gamma})$. Let $t\in(0,T]$. Using \eqref{Est-MLF} we get that there is a constant $C>0$ such that,
\begin{align}\label{Est-1}
\norm{S_1(t)u_0}^2_{V_\gamma}\le 4\sum_{n=1}^\infty \,\abs{u_{0,n}\lambda_n^{\gamma} t^{-\beta}E_{\mu,1-\beta}(-\lambda_nt^\mu)}^2 \le Ct^{-2\beta}\norm{u_0}_{V_\gamma}^2.
\end{align}
Using \eqref{Est-ML2} we obtain that there is a constant $C>0$ such that for every $t\in(0,T]$ (recall that $\mu\gamma=1$),
\begin{align}\label{Est-2}
\norm{S_2(t)u_1}^2_{V_\gamma}\le 4\sum_{n=1}^\infty\,\abs{u_{1,n}\lambda^{\gamma}t^{1-\beta}E_{\mu,2-\beta}(-\lambda_nt^\mu)}^2\le Ct^{-2\beta}\norm{u_1}_{L^2(\Omega)}^2.
\end{align}
Using \eqref{Est-ML2} again,  the Minkowski and the H\"older inequalities, we obtain that there is a constant $C>0$ such that for every $t\in(0,T]$,
\begin{align}\label{Est-3}
\norm{S_3(t)f}_{V_\gamma}&\le 2\,\int_0^t\left( \sum_{n=1}^\infty\abs{f_n(\tau)\lambda_n^\gamma(t-\tau)^{\mu-1}E_{\mu,\mu}(-\lambda_n(t-\tau)^\mu)}^2\right)^\frac 12 d\tau\notag \\&\le C\int_0^t\,(t-\tau)^{\mu-2}\left( \sum_{n=1}^\infty\abs{f_n(\tau)}^2\right) ^\frac 12 d\tau \notag\\&\le C\int_0^t(t-\tau)^{\mu-2}\norm{f(\cdot,\tau)}_{L^2(\Omega)}d\tau \notag\\&\le Ct^{\mu-1-\frac 1p}\norm{f}_{L^p((0,T);L^2(\Omega))}.
\end{align}
 It follows from  \eqref{Est-1}, \eqref{Est-2} and \eqref{Est-3} that there is a constant $C_1>0$ such that for every $t\in (0,T]$, 
\begin{align}\label{CDS}
\norm{u(\cdot,t)}_{V_\gamma}\le C_1\left( t^{-\beta}\norm{u_0}_{V_\gamma}+t^{-\beta}\norm{u_1}_{L^2(\Omega)}+t^{\mu-1-\frac 1p}\norm{f}_{L^p((0,T);L^2(\Omega))}\right).
\end{align}
It is straightforward to show that the series \eqref{Weak-Sol} converges in $V_\gamma$ uniformly  in compact subsets of $(0,T]$.  Hence,  $u\in C((0,T];V_\gamma)$.

{\bf Step 3}:  Next, we show that $I_t^{\beta}u\in C([0,T];V_\gamma)$. Using Lemma \ref{int-frac} and the identity \eqref{Con-Int} we have that
\begin{align}\label{Init-Cond1}
I_t^\beta u(\cdot,t)=&\sum_{n=1}^\infty u_{0,n}E_{\mu,1}(-\lambda_nt^\mu)\varphi_n+\sum_{n=1}^\infty u_{1,n}tE_{\mu,2}(-\lambda_nt^\mu)\varphi_n\notag\\&\;\;\;+\sum_{n=1}^\infty\left( \int_0^t f_n(\tau)(t-\tau)^{1-\nu(2-\mu)}E_{\mu,2-\nu(2-\mu)}(-\lambda_n(t-\tau)^\mu)\,d\tau\right) \varphi_n.
\end{align}
Proceeding as in Step 2, we get that there is a constant $C>0$ such that  the following estimates hold:
\begin{align*}
\norm{I_t^\beta S_1(t)u_0}_{V_\gamma}\le C\norm{u_0}_{V_\gamma} \;\mbox{and}\; \norm{I_t^\beta S_2(t)u_1}_{V_\gamma}\le C\norm{u_1}_{L^2(\Omega)}, \;\forall t\in[0,T]  
\end{align*}
and
\begin{align}\label{Est-4}
\norm{I_t^\beta S_3(t)f}_{V_\gamma}\le Ct^{\nu(\mu-2)+1-\frac 1p}\norm{f}_{L^p((0,T);L^2(\Omega))}, \;\forall t\in[0,T].
\end{align}
Here also, it is straightforward to show that the series \eqref{Init-Cond1} converges in $V_\gamma$ uniformly in $t\in[0,T]$. Hence,  $I_t^\beta u\in C([0,T];V_\gamma)$.

{\bf Step 4}:  Next, we show that $I_t^\beta u\in C^1([0,T];L^2(\Omega))$.  Using Step 3 and the continuous embedding $V_\gamma\hookrightarrow L^2(\Omega)$,  it suffices to show that $\partial_t I_t^\beta u\in C([0,T];L^2(\Omega))$.
A simple calculation gives for $t\in[0,T]$,
\begin{align}\label{Init-Cond2}
\partial_tI_t^\beta u(\cdot,t)=&\sum_{n=1}^\infty u_{0,n}\lambda_nt^{\mu-1}E_{\mu,\mu}(-\lambda_nt^\mu)\varphi_n+\sum_{n=1}^\infty u_{1,n}E_{\mu,1}(-\lambda_nt^\mu)\varphi_n\notag\\
&+\sum_{n=1}^\infty\left( \int_0^t\,f_n(\tau)(t-\tau)^{\nu(\mu-2)}E_{\mu,\nu(\mu-2)+1}(-\lambda_n(t-\tau)^\mu)d\tau\right) \varphi_n\notag\\&=\partial_tI^\beta_tS_1(t)u_0+\partial_tI^\beta_tS_2(t)u_1+\partial_tI^\beta_t(S_3\ast f)(t).
\end{align}
Using \eqref{Est-ML2}, \eqref{Est-ML3} and proceeding as in Step 3  we obtain that there is a constant $C>0$ such that  the following estimates hold:
\begin{align}\label{MJW-K}
\norm{\partial_tI^\beta_tS_1(t)u_0}_{L^2(\Omega)}&\le C\norm{u_0}_{V_\gamma}\;\;\;\mbox{and}\;\;\;\norm{\partial_tI^\beta_tS_2(t)u_1}_{L^2(\Omega)}\le C\norm{u_1}_{L^2(\Omega)},\;\; \forall t\in[0,T],
\end{align} 
and
\begin{align}\label{Est-5}
\norm{\partial_tI^\beta_t(S_3\ast f)(t)}_{L^2(\Omega)}\le Ct^{\nu(\mu-2)+1-\frac 1p}\norm{f}_{L^p((0,T);L^2(\Omega))},\;\; \forall t\in[0,T].
\end{align}
Notice that, since  for $\nu\neq 1$ and $\mu\neq 2$, $\nu(\mu-2)+2>\mu$ and by assumption $p>1/(\mu-1)$, we have that $\nu(\mu-2)+1-1/p>\mu-1-1/p>0$.
Since the series \eqref{Init-Cond2} converges in $L^2(\Omega)$ uniformly in $t\in[0,T]$,  we have that $I_t^\beta u\in C^1([0,T];L^2(\Omega))$.  It follows from  \eqref{Est-4},  \eqref{MJW-K}  and \eqref{Est-5} that there is a constant $C>0$ such that
\begin{align*}
\norm{I_t^\beta u(\cdot,t)}_{V_\gamma}+\norm{\partial_tI_t^\beta u(\cdot,t)}_{L^2(\Omega)}\le C\left( \norm{u_0}_{V_\gamma}+\norm{u_1}_{L^2(\Omega)}+t^{\nu(\mu-2)+1-\frac 1p}\norm{f}_{L^p((0,T);L^2(\Omega))}\right) .
\end{align*}

{\bf Step 5}:  We show that $\D_t^{\mu,\nu}u\in C((0,T];V_{-\gamma})$. It follows from \eqref{Weak-Sol} that 
\begin{align}\label{sd}
\mathbb{D}^{\mu,\nu}u(\cdot,t)=&-\sum_{n=1}^\infty u_{0,n}\lambda_nt^{-\beta}E_{\mu,1-\beta}(-\lambda_nt^\mu)\varphi_n-\sum_{n=1}^\infty u_{1,n}\lambda_nt^{1-\beta}E_{\mu,2-\beta}(-\lambda_nt^\mu)\varphi_n \notag\\ &-\sum_{n=1}^\infty \left( \int_0^t (f(\cdot, \tau), \varphi_n)\lambda_n(t-\tau)^{\mu-1}E_{\mu,\mu}(-\lambda_n(t-\tau)^\mu)d\tau \right) \varphi_n+f(\cdot,t).
\end{align}
Using \eqref{Est-MLF},  \eqref{Est-ML3} and proceeding as above,  we get that there is  constant $C>0$ such that  the following estimates hold (recall that $\lambda_1>0$ and $1<2\gamma=\frac 2\mu<2$):
\begin{align}
\norm{\sum_{n=1}^\infty u_{0,n}\lambda_nt^{-\beta}E_{\mu,1-\beta}(-\lambda_nt^\mu)\varphi_n}_{V_{-\gamma}}^2\le& C\sum_{n=1}^\infty |u_{0,n}\lambda_n^{1-\gamma}t^{-\beta}E_{\mu,1-\beta}(-\lambda_nt^\mu)|^2\notag\\
\le &\sum_{n=1}^\infty |u_{0,n}\lambda_n^\gamma\lambda_n^{1-2\gamma}t^{-\beta}E_{\mu,1-\beta}(-\lambda_nt^\mu)|^2\notag\\
\le &C\lambda_1^{1-2\gamma}t^{-2\beta}\norm{u_0}_{V_\gamma}^2,
\end{align}
and
\begin{align}
\norm{\sum_{n=1}^\infty u_{1,n}\lambda_nt^{1-\beta}E_{\mu,2-\beta}(-\lambda_nt^\mu)\varphi_n}_{V_{-\gamma}}^2\le &C\sum_{n=1}^\infty |u_{1,n}\lambda_n^{1-\gamma}t^{1-\beta}E_{\mu,2-\beta}(-\lambda_nt^\mu)|^2\notag\\
\le &Ct^{2\nu(2-\mu)}\norm{u_1}_{L^2(\Omega)}^2.
\end{align}
Similarly,  using  \eqref{Est-ML3} we can easily deduce that
\begin{align}
\norm{\sum_{n=1}^\infty \left( \int_0^t (f(\cdot, \tau), \varphi_n)\lambda_n (t-\tau)^{\mu-1}E_{\mu,\mu}(-\lambda_n(t-\tau)^\mu)d\tau \right) \varphi_n}_{V_{-\gamma}}\le Ct^{1-\frac 1p}\norm{f}_{L^p((0,T);L^2(\Omega))}.
\end{align}
Since the series \eqref{sd} converges in $V_{-\gamma}$ uniformly in compact subsets of $(0,T]$,  we can conclude that $\D_t^{\mu,\nu}u\in C((0,T];V_{-\gamma})$.

{\bf Step 6}:  Since $\D_t^{\mu,\nu}u\in C((0,T];V_{-\gamma})$, $Au(\cdot,t)\in V_{-1/2}\subset V_{-\gamma}$ and $f(\cdot,t)\in L^2(\Omega)$ for a.e. $t\in(0,T)$, then taking the duality product in \eqref{sd} with $\varphi\in V_{\gamma}$, we get the variational identity \eqref{Var-I}.

{\bf Step 7}:  Using \eqref{Init-Cond1} and \eqref{Init-Cond2}, we obtain that
\begin{equation}
I^\beta_tu(\cdot,0)=\sum_{n=1}^\infty u_{0,n}\varphi_n=u_0\;\;\mbox{and}\;\;\partial_tI^\beta_tu(\cdot,0)=\sum_{n=1}^\infty u_{0,n}\varphi_n=u_1,
\end{equation}
and we have shown \eqref{Init-Cond}. 

{\bf Step 8}:  Finally,  we show uniqueness. Let $u,\,v$ be two weak solutions of \eqref{main-EQ} with the same initial data $u_0$, $u_1$ and source $f$.  We can deduce from \eqref{CDS} (the continuous dependence of solutions on the given data) that $u=v$ and the proof is finished.
\end{proof}

%
%

\subsection{Another representation of solutions to the homogeneous equation}
Let $0\le\nu\le$, $1<\mu\le 2$ and consider the following homogeneous equation:
\begin{equation}\label{main-EQ1}
\begin{cases}
\mathbb D_t^{\mu,\nu}u +Au=0\;\;&\mbox{ in }\; \Omega\times (0,T),\\
(I_t^{(1-\nu)(2-\mu)}u)(\cdot,0)=u_0 &\mbox{ in }\;\Omega,\\
\left( \partial_t I_t^{(1-\nu)(2-\mu)}u\right)(\cdot,0)=u_1 &\mbox{ in }\;\Omega.
\end{cases}
\end{equation}
The main concern of this section is to give another useful representation  (apart from  \eqref{Weak-Sol}) of solutions to \eqref{main-EQ1}.  For this we define the following operators.

\begin{definition} 
Let $1<\mu\le 2$.
Given $u\in L^2(\Omega)$ and $t\ge 0$, we let
\begin{align}\label{Ope}
S_\mu(t)u:=\sum_{n=1}^{\infty}(u,\varphi_n)_{L^2(\Om)}E_{\mu,\mu}(-\lambda_nt^\mu)\varphi_n
\end{align}
and
\begin{align}\label{Ope-2}
S_{\mu-1}(t)u:=\sum_{n=1}^\infty\,(u,\varphi_n)_{L^2(\Om)}E_{\mu,\mu-1}(-\lambda_nt^\mu)\varphi_n.
\end{align}
\end{definition}

We have the following results.

\begin{lemma}\label{family-op} 
Let  $S_\mu(t)$ be the operator defined in \eqref{Ope}. Then the following assertions hold.
\begin{enumerate}
\item There is a constant $C_1>0$ such that for every $u\in L^2(\Omega)$ and $t\ge 0$, 
\begin{align*}
\norm{S_\mu(t)u}_{L^2(\Omega)}\le C_1\norm{u}_{L^2(\Omega)}.
\end{align*}
\item For every $u\in D(A)$, we have that $S_\mu(t)u\in D(A)$  and $AS_\mu(t)u=S_\mu(t)Au$  for all $t\ge 0$.
\item $S_\mu(t)S_\mu(\tau)=S_\mu(\tau)S_\mu(t)$ for all $t, \tau\geq0$.
\item There is a constant $C_2>0$ such that for every $u\in L^2(\Omega)$ and $t>0$, 
$$\norm{\dfrac{dS_\mu(t)u}{dt}}_{L^2(\Om)}\le C_2t^{-1}\norm{u}_{L^2(\Omega)}.$$
\end{enumerate}
\end{lemma}

\begin{proof}
(a) This assertion follows directly from the definition of $S_\mu(t)$  in \eqref{Ope} and the estimate \eqref{Est-MLF} of the Mittag-Leffler function.

(b) Let $u\in D(A)$ and $t\ge 0$. Using \eqref{Est-MLF}, we get that there is a constant $C>0$ such that for every $t\ge 0$, 
\begin{align*}
\norm{S_\mu(t)u}^2_{D(A)}\le&\sum_{n=1}^\infty\abs{\lambda_n(S_\mu(t)u,\varphi_n)_{L^2(\Om)}}^2=\sum_{n=1}^\infty\abs{\lambda_n(u,\varphi_n)_{L^2(\Om)}E_{\mu,\mu}(-\lambda_nt^\alpha)}^2\le C\norm{u}^2_{D(A)}. 
\end{align*}
Thus, $S_\mu(t)u\in D(A)$.  The second part of the assertion is easy to see.

(c) This part is obtained by a simple calculation and using that $(\varphi_n)_{n\in \mathbb{N}}$ is an orthonormal basis of $L^2(\Om)$.

(d) Let $u\in L^2(\Omega)$ and $t>0$. Since the series  \eqref{Ope} converges in $L^2(\Om)$ uniformly in compact subsets of $[0,\infty)$,  we have that
\begin{align}\label{DS-1}
\dfrac{dS_\mu(t)u}{dt}=\sum_{n=1}^\infty(u,\varphi_n)_{L^2(\Om)}\dfrac{d}{dt}\Big(E_{\mu,\mu}(-\lambda_nt^\mu)\Big)\varphi_n.
\end{align}
Using \eqref{MLD} and the Chain-rule,  we obtain that 
 \begin{align*}
\dfrac{d}{dt}\Big[E_{\mu,\mu}(-\lambda_nt^\mu)\Big]=\frac{E_{\mu,\mu-1}(-\lambda_nt^\mu)+(1-\mu)E_{\mu,\mu}(-\lambda_nt^\mu)}{-\mu\lambda_nt^\mu} (-\lambda_n\mu t^{\mu-1}).
\end{align*}
Using  \eqref{Est-MLF}  and the previous identity we have that there is a constant $C>0$ such that 
\begin{align}\label{DS-2}
\abs{\dfrac{d}{dt}\Big[E_{\mu,\mu}(-\lambda_nt^\mu)\Big]}\le Ct^{-1}.
\end{align}
Thus, the assertion follows by combining \eqref{DS-1}-\eqref{DS-2}. The proof is finished.
\end{proof}

\begin{lemma}\label{family-op-2} 
Let  $S_{\mu-1}(t)$ be the operator defined in \eqref{Ope-2}.
Then,  the following assertions hold.
\begin{enumerate}
\item There is a constant $C_1>0$ such that for every $u\in L^2(\Omega)$ and $t\ge 0$, 
\begin{align}
\norm{S_{\mu-1}(t)u}_{L^2(\Omega)}\le C_1\norm{u}_{L^2(\Omega)}.
\end{align}

\item For every $u\in D(A)$, we have that $S_{\mu-1}(t)u\in D(A)$ and $AS_{\mu-1}(t)u=S_{\mu-1}(t)Au$ for all $t\ge 0$.

\item $S_{\mu-1}(t)S_{\mu-1}(\tau)=S_{\mu-1}(\tau)S_{\mu-1}(t)$ for all $t, \tau\geq0$.

\item There is a constant $C_2>0$ such that for every $u\in L^2(\Omega)$ and $t>0$, 
$$\norm{\dfrac{dS_{\mu-1}(t)u}{dt}}_{L^2(\Om)}\le C_2t^{-1}\norm{u}_{L^2(\Omega)}.$$
\end{enumerate}
\end{lemma}

\begin{proof}
The proof follows as the proof of Lemma \ref{family-op}.  For Part (d) one has to use the Chain-rule and the fact that
\begin{align*}
\frac{d}{dz}E_{\mu,\mu-1}(z)=\frac{E_{\mu,\mu-2}(z)-(\mu-2)E_{\mu,\mu-1}(z)}{\mu z}.
\end{align*}
We omit the details.
\end{proof}

We have the following representation of weak solutions of the homogeneous equation.

\begin{theorem}\label{classical}
Let $0\le\nu\le 1$, $1<\mu\le 2$, $\gamma=1/\mu$ or $\gamma=1/2$, $u_0\in V_\gamma$, $u_1\in L^2(\Omega)$ and let $u$ be the unique weak solution of  \eqref{main-EQ1} which is given by
\begin{align}\label{class-sol}
u(\cdot,t):=&\sum_{n=1}^\infty(u_0,\varphi_n)_{L^2(\Om)}t^{-(1-\nu)(2-\mu)}E_{\mu,\mu-1+\nu(2-\mu)}(-\lambda_nt^\mu)\varphi_n\notag\\
&+\sum_{n=1}^\infty(u_1,\varphi_n)_{L^2(\Om)}t^{\mu-1+\nu(2-\mu)}E_{\mu,\mu+\nu(2-\mu)}(-\lambda_nt^\mu)\varphi_n.
\end{align}
Then
\begin{align}\label{e248}
u(\cdot,t):=I_t^{\nu(2-\mu)}t^{\mu-2}S_{\mu-1}(t)u_0+I_t^{\nu(2-\mu)}t^{\mu-1}S_\mu(t)u_1.
\end{align}
\end{theorem}

\begin{proof}
The representation \eqref{e248} follows from \eqref{Ope}, \eqref{Ope-2}, \eqref{class-sol} by using the method of Laplace transform.  We omit the details.
\end{proof}

\subsection{Well-posedness of the adjoint system}
This section is devoted to the existence, regularity and uniqueness of solutions to the dual system \eqref{ACP-Dual}.
Observe that in  \eqref{ACP-Dual} the final conditions have the interpretation:
\begin{align*}
\left(I_{t,T}^{\nu(2-\mu)}v\right)(\cdot,T)=\lim_{t\to T}\frac{1}{\Gamma(\nu(2-\mu))}\,\int_t^T\,(\tau-t)^{\nu(2-\mu)-1}v(\cdot,\tau)\,d\tau=v_0,
\end{align*}
and
\begin{align*}
\left(D_{t,T}^{1-\nu(2-\mu)}v\right)(\cdot,T)=-\lim_{t\to T}\frac{1}{\Gamma(\nu(2-\mu))}\dfrac{d}{dt}\,\int_t^T\,(\tau-t)^{\nu(2-\mu)-1}v(\cdot,\tau)\,d\tau=v_1,
\end{align*}

We adopt the following definition of weak solutions.

\begin{definition}
Let $0\le  \nu\le 1$, $1<\mu \le 2$, $1/2\le\gamma\le 1$, $v_0\in V_\gamma$ and $v_1\in L^2(\Omega)$.
A function $v$ is said to be a weak solution of \eqref{ACP-Dual}, if  the following properties hold:
\begin{itemize}
\item Regularity:
\begin{equation}\label{Dual-egu}
\begin{cases}
v\in C([0,T);V_\gamma),\\
 I_{t,T}^{\nu(2-\mu)}v\in C([0,T];V_\gamma),\\
 D_{t,T}^{1-\nu(2-\mu)}v\in C([0,T];L^2(\Omega)), \\
\mathbb D_{t,T}^{\mu,1-\nu} v\in C([0,T);V_{-\gamma}).
\end{cases}
\end{equation}
\item Final conditions:
\begin{align*}
\left(  I_{t,T}^{\nu(2-\mu)}v\right)(\cdot,T)=v_0 \;\mbox{ and } \; \left(D_{t,T}^{1-\nu(2-\mu)}v\right)(\cdot,T)=v_1 \;\mbox{ a.e.  in } \Omega.
\end{align*}
\item Variational  identity: For every $\varphi\in V_\gamma$ and a.e. $t\in (0,T)$, 
\begin{align}\label{Dual-Var-I}
\langle \mathbb D_{t,T}^{\mu,(1-\nu)} v(\cdot,t),\varphi\rangle_{V_{-\gamma},V_\gamma} +\mathcal E_A(v(\cdot,t),\varphi)=0.
\end{align}
\end{itemize}
\end{definition}

We have the following existence and uniqueness result.

\begin{theorem}\label{Dual-theo-weak}
Let $0\le  \nu\le 1$, $1<\mu \le 2$, $\gamma=1/\mu$ or $\gamma=1/2$, $v_0\in V_\gamma$ and $v_1\in L^2(\Omega)$. Then the  system \eqref{ACP-Dual} has  a unique weak solution $v$ given by
 \begin{align}\label{Dual-sol-spec}
 v(\cdot,t)=&\sum_{n=1}^\infty (v_0,\varphi_n)_{L^2(\Om)}(T-t)^{\nu(\mu-2)}E_{\mu,\nu(\mu-2)+1}(-\lambda_n(T-t)^\mu)\varphi_n \notag\\ 
 &+\sum_{n=1}^\infty (v_1,\varphi_n)_{L^2(\Om)}(T-t)^{\nu(\mu-2)+1}E_{\mu,\nu(\mu-2)+2}(-\lambda_n(T-t)^\mu)\varphi_n.
 \end{align}
  Moreover,  $v$ can be analytically extended to the half-plane 
 $$\Sigma_T:=\{z\in\CC:\;\mbox{Re}(z)<T\}.$$
\end{theorem}

\begin{proof}
We proceed in several steps.

{\bf Step 1}: First, we show uniqueness of solutions. Indeed, let $v$ be a solution of  \eqref{ACP-Dual} with $v_0=v_1=0$. Taking the duality product of \eqref{ACP-Dual} with $\varphi_n$ and setting $v_n(t):=(v(t),\varphi_n)_{L^2(\Omega)}$, we obtain that (given that $A$ is a selfadjoint operator) for every $n\in\mathbb N$,
\begin{equation}\label{O1}
\mathbb D_{t,T}^{\mu, 1-\nu} v_n(t)=-\lambda_nv_n(t),\;\;\mbox{ for a.e. }\; t\in (0,T).
\end{equation}
Since $I_{t,T}^{\nu(2-\mu)}v\in C([0,T];L^2(\Omega))$, we have that $I_{t,T}^{\nu(2-\mu)}v_n=(I_{t,T}^{\nu(2-\mu)}v,\varphi_n)_{L^2(\Omega)}\in C[0,T]$ and
\begin{align*}
|I_{t,T}^{\nu(2-\mu)}v_n(t)|^2\le \sum_{n=1}^\infty|I_{t,T}^{\nu(2-\mu)}v_n|^2\le \|I_{t,T}^{\nu(2-\mu)}v\|_{L^2(\Omega)}^2\longrightarrow 0\;\mbox{ as }\; t\to T.
\end{align*}
This implies that
\begin{equation}\label{O2}
I_{t,T}^{\nu(2-\mu)}v_n(T)=0.
\end{equation}
 Similarly we have that
\begin{equation} \label{O3}
 D_{t,T}^{1-\nu(2-\mu)}v_n(T)=0.
 \end{equation}
Since the fractional ordinary differential equation \eqref{O1} with the final conditions \eqref{O2}-\eqref{O3} has a unique solution $v_n$ given by (this can be easily shown by using the method of Laplace transform)
\begin{align*}
v_n(t)=&(T-t)^{\nu(\mu-2)}E_{\mu,\nu(\mu-2)+1}(-\lambda_n(T-t)^\mu)I_{t,T}^{\nu(2-\mu)}v_n(T)\\&\;\;\;\;\;\;\;+(T-t)^{\nu(\mu-2)+1}E_{\mu,\nu(\mu-2)+2}(-\lambda_n(T-t)^\mu)D_{t,T}^{1-\nu(2-\mu)}v_n(T),
\end{align*}
it follows that $v_n(t)=0$ for every $n\in\NN$. Since $(\varphi_n)$ is a complete orthonormal system in $L^2(\Omega)$, we have that $v=0$ in $\Omega\times (0,T)$ and the proof of the uniqueness is complete.

{\bf Step 2}: Next, we prove the existence of solutions. Let  $v_{0,n}:=(v_0,\varphi_n)_{L^2(\Omega)}$, $v_{1,n}:=(v_1,\varphi_n)_{L^2(\Omega)}$, $1\le n\le k$ where $n,k\in\NN$,  and set
\begin{align}\label{vk}
v_k(x,t):=&\sum_{n=1}^k v_{0,n}(T-t)^{\nu(\mu-2)}E_{\mu,\nu(\mu-2)+1}(-\lambda_n(T-t)^\mu)\varphi_n(x)\notag\\
 &+\sum_{n=1}^kv_{1,n}(T-t)^{\nu(\mu-2)+1}E_{\mu,\nu(\mu-2)+2}(-\lambda_n(T-t)^\mu)\varphi_n(x).
\end{align}
\begin{enumerate}
\item Let $v$ be given by \eqref{Dual-sol-spec}. We claim that $v\in C([0,T);V_\gamma)$. Using \eqref{Est-MLF} and the estimates in Lemma \ref{lem-INE}, we have that for every $t\in [0,T)$,
\begin{align*}
\|v_k(\cdot,t)-v_m(\cdot,t)\|_{V_\gamma}^2=&2\sum_{n=k+1}^m\left| \lambda^\gamma_n v_{0,n}(T-t)^{\nu(\mu-2)}E_{\mu,\nu(\mu-2)+1}(-\lambda_n(T-t)^\mu)\right|^2\\
&+2\sum_{n=k+1}^m\left| v_{1,n}\lambda_n^\gamma(T-t)^{\nu(\mu-2)+1}E_{\mu,\nu(\mu-2)+2}(-\lambda_n(T-t)^\mu)\right|^2\\
&\le (T-t)^{2\nu(\mu-2)}\sum_{n=k+1}^m\left| \lambda^\gamma_n v_{0,n}\right|^2 +(T-t)^{\nu(\mu-2)}\sum_{n=k+1}^m\left| v_{1,n}\right|^2\\
&\longrightarrow 0\;\mbox{ as }\; k,m\to\infty.
\end{align*}
We have shown that the series \eqref{Dual-sol-spec} converges in $V_\gamma$ and the convergence is uniform in compact subsets of $[0,T)$.  Hence, $v\in C([0,T);V_\gamma)$. In addition,  there is a constant $C>0$ such that for every $t\in [0,T)$,
\begin{align}\label{norm-V}
\|v(\cdot,t)\|_{V_\gamma}^2\le C(T-t)^{2\nu(\mu-2)}\left(\|v_0\|_{\gamma}^2+\|v_1\|_{L^2(\Omega)}^2\right).
\end{align}
\item Let $v$ be given by \eqref{Dual-sol-spec}. We claim that  $I_{t,T}^{\nu(2-\mu)}v\in C([0,T];V_\gamma)$. Integrating \eqref{vk} termwise and using \eqref{EE-11}, we get
\begin{align}\label{V1}
I_{t,T}^{\nu(2-\mu)}v_k(\cdot,t)=&\sum_{n=1}^k v_{0,n}E_{\mu,1}(-\lambda_n(T-t)^\mu)\varphi_n\notag\\
& +\sum_{n=1}^kv_{1,n}\int_t^TE_{\mu,1}(-\lambda_n(T-\tau)^\mu)\;d\tau\varphi_n.
\end{align}
Using \eqref{Est-MLF} and the estimates in Lemma \ref{lem-INE}, we have that for every $t\in [0,T]$ and $m,k\in\NN$ with $m>k$,
\begin{align*}
\|I_{t,T}^{2-\alpha}v_k(\cdot,t)-I_{t,T}^{\nu(2-\mu)}v_m(\cdot,t)\|_{V_\gamma}^2=&2\sum_{n=k+1}^m |\lambda_n^\gamma
v_{0,n}E_{\mu,1}(-\lambda_n(T-t)^\mu)|^2\\
& +2\sum_{n=k+1}^m|v_{1,n}\lambda_n^{\gamma}\int_t^TE_{\mu,1}(-\lambda_n(T-\tau)^\mu)\;d\tau|^2\\
=&2\sum_{n=k+1}^m |\lambda_n^\gamma
v_{0,n}E_{\mu,1}(-\lambda_n(T-t)^\mu)|^2\\
& +2\sum_{n=k+1}^m|v_{1,n}\lambda_n^{\gamma}(T-t)E_{\mu,2}(-\lambda_n(T-t)^\mu)|^2\\
\le &C\left(\sum_{n=k+1}^m |\lambda_n^\gamma v_{0,n}|^2+\sum_{n=k+1}^m|v_{1,n}|^2\right)\\
&\longrightarrow 0\;\mbox{ as }\; k,m\to\infty.
\end{align*}
We have shown that the series
\begin{align*}
&\sum_{n=1}^\infty v_{0,n}E_{\mu,1}(-\lambda_n(T-t)^\mu)\varphi_n\notag\\
 +&\sum_{n=1}^\infty v_{1,n}\int_t^TE_{\mu,1}(-\lambda_n(T-\tau)^\mu)\;d\tau\varphi_n \longrightarrow I_{t,T}^{\nu(2-\mu)}v(\cdot,t)\;\mbox{ in }\;V_\gamma,
\end{align*}
and that the convergence is uniform in $t\in [0,T]$. Hence, $I_{t,T}^{\nu(2-\mu)}v\in C([0,T];V_\gamma)$. Using \eqref{Est-MLF} and Lemma \ref{lem-INE} again, we obtain that there is a constant $C>0$ such that for every $t\in [0,T]$,
\begin{align*}
 \|I_{t,T}^{\nu(2-\mu)}v(\cdot,t)\|_{V_{\gamma}}^2\le C^2(\|v_0\|_{V_\gamma}^2+\|v_1\|_{L^2(\Omega)}^2).
\end{align*}

\item We show that $D_{t,T}^{1-\nu(2-\mu)}v\in C([0,T];L^2(\Omega))$. Integrating \eqref{vk} termwise and using Lemma \ref{Der-MLF}, we also get 
\begin{align}\label{V2}
D_{t,T}^{1-\nu(2-\mu)}v_k(t)=&\sum_{n=1}^kv_{0,n}\lambda_n(T-t)^{\mu -1}E_{\mu,\mu}(-\lambda_n(T-t)^\mu)\varphi_n\notag\\
& +\sum_{n=1}^kv_{1,n}E_{\mu,1}(-\lambda_n(T-t)^\mu)\varphi_n.
\end{align}
Proceeding as in Part (a) or (b),  we obtain that 
\begin{align*}
&\sum_{n=1}^\infty v_{0,n}\lambda_n(T-t)^{\mu -1}E_{\mu,\mu}(-\lambda_n(T-t)^\mu)\varphi_n\\
 +&\sum_{n=1}^\infty v_{1,n}E_{\mu,1}(-\lambda_n(T-t)^\mu)\varphi_n\longrightarrow D_{t,T}^{1-\nu(2-\mu)}v(t)\;\mbox{ in }\;L^2(\Omega),
\end{align*}
and the convergence is uniform in $t\in [0,T]$. Hence, $D_{t,T}^{1-\nu(2-\mu)}v\in C([0,T];L^2(\Om))$.
Using \eqref{Est-MLF} and  Lemma \ref{lem-INE} we obtain that there is a  constant $C>0$ such that for every $t\in[0,T]$,
\begin{align}\label{254}
\|D_{t,T}^{1-\nu(2-\mu)}v(\cdot,t)\|_{L^2(\Omega)}^2\le &2\sum_{n=1}^\infty |\lambda_n^{\gamma}v_{0,n}|^2|\lambda_n^{1-\gamma}(T-t)^{\mu -1}E_{\mu,\mu}(-\lambda_n(T-t)^\mu)|^2\notag\\
& +2\sum_{n=1}^\infty |v_{1,n}|^2|E_{\mu,1}(-\lambda_n(T-t)^\mu)|^2\notag\\
\le &C\left(\|v_0\|_{V_\gamma}^2+C_2^2\|v_1\|_{L^2(\Omega)}^2\right).
\end{align}

\item We prove that $\D_{t,T}^{\mu,1-\nu} v\in C([0,T);L^2(\Omega))$.  Proceeding as above we can deduce that $\mathbb{D}_{t,T}^{\mu,1-\nu}v\in C([0,T);L^2(\Omega)$.
Since $\D_{t,T}^{\mu,1-\nu} v(\cdot,t)=-Av(\cdot,t)$, applying  Lemma \ref{lem-INE} again,  we have that there is a constant $C>0$ such that for every $t\in [0,T)$,
\begin{align*}
\|\D_{t,T}^{\mu,1-\nu}v(\cdot,t)\|_{L^2(\Omega)}^2\le &\sum_{n=1}^\infty |\lambda_n^{\gamma}v_{0,n}|^2|\lambda_n^{1-\gamma}(T-t)^{\nu(\mu -2)}E_{\mu,\nu(\mu-2)+1}(-\lambda_n(T-t)^\mu)|^2\notag\\
& +\sum_{n=1}^\infty |v_{1,n}|^2|\lambda_n(T-t)^{\nu(\mu-2)+1}E_{\mu,\nu(\mu-2)+2}(-\lambda_n(T-t)^\mu)|^2\notag\\
\le &C(T-t)^{-2+2(1-\nu)(2-\mu)}\left(\|v_0\|_{V_\gamma}^2+\|v_1\|_{L^2(\Omega)}^2\right).
\end{align*}

\item It follows from \eqref{V1} and \eqref{V2} that
\begin{align*}
I_{t,T}^{\nu(2-\mu)}v(\cdot,T)=v_0\;\mbox{ and }\; D_{t,T}^{1-\nu(2-\mu)}v(\cdot,T)=v_1\;\mbox{ a.e. in }\;\Omega.
\end{align*}
The proof of the existence is complete.
\end{enumerate}

{\bf Step 3}: Finally, since $E_{\mu,\nu(\mu-2)+1}(-\lambda_nz)$ and $E_{\mu,\nu(\mu-2)+2}(-\lambda_nz)$ are entire functions, it follows that the two functions
\begin{align*}
(T-t)^{\nu(\mu-2)}E_{\mu,\nu(\mu-2)+1}(-\lambda_n(T-t)^\mu)\;\mbox{ and }\;
(T-t)^{\nu(\mu-2)+1}E_{\mu,\nu(\mu-2)+2}(-\lambda_n(T-t)^\mu)
\end{align*}
 can be analytically extended to the half-plane $\Sigma_T$. This implies that the functions
 \begin{align*}
 \sum_{n=1}^k v_{0,n}(T-t)^{\nu(\mu-2)}E_{\mu,\nu(\mu-2)+1}(-\lambda_n(T-t)^\mu)\varphi_n
\end{align*}
 and
\begin{align*}
 \sum_{n=1}^kv_{1,n}(T-t)^{\nu(\mu-2)+1}E_{\mu,\nu(\mu-2)+2}(-\lambda_n(T-t)^\mu)\varphi_n
\end{align*}
  are analytic in $\Sigma_T$. Let $\delta>0$ be fixed but otherwise arbitrary. Let $z\in\CC$ satisfy $\mbox{Re}(z)\le T-\delta$. Then using Lemma \ref{lem-INE}, we obtain that there is a constant $C>0$ such that 
 \begin{align*}
& \left\Vert\sum_{n=k+1}^\infty v_{0,n}(T-z)^{\nu(\mu-2)}E_{\mu,\nu(\mu-2)+1}(-\lambda_n(T-z)^\mu)\varphi_n\right\Vert_{L^2(\Omega)}^2\\
 &+\left\Vert\sum_{n=k+1}^\infty v_{1,n}(T-z)^{\nu(\mu-2)+1}E_{\mu,\nu(\mu-2)+2}(-\lambda_n(T-z)^\mu)\varphi_n\right\Vert_{L^2(\Omega)}^2\\
 \le &C\sum_{n=k+1}^\infty |v_{0,n}|^2|T-z|^{2\nu(\mu-2)}\left(\frac{1}{1+\lambda_n|T-z|^{\mu}}\right)^2\\
 &+C\sum_{n=k+1}^\infty |v_{1,n}|^2|T-z|^{2(\nu(\mu-2)+1)}\left(\frac{1}{1+\lambda_n|T-z|^{\mu}}\right)^2\\
 \le &C\delta^{2\nu(\mu-2)}\sum_{n=k+1}^\infty |v_{0,n}|^2 +C\delta^{2(\nu(\mu-2)+1)}\sum_{n=k+1}^\infty |v_{1,n}|^2\longrightarrow 0\;\mbox{ as }\; k\to\infty.
 \end{align*}
 We have shown that
 \begin{align*}
 v(\cdot,z):=&\sum_{n=1}^\infty (v_0,\varphi_n)(T-z)^{\nu(\mu-2)}E_{\mu,\nu(\mu-2)+1}(-\lambda_n(T-z)^\mu)\varphi_n\notag\\
 &+\sum_{n=1}^\infty(v_1,\varphi_n)(T-z)^{\nu(\mu-2)+1}E_{\mu,\nu(\mu-2)+2}(-\lambda_n(T-z)^\mu)\varphi_n
 \end{align*}
 is uniformly convergent in any compact subset of $\Sigma_T$. Hence, $v$ is also analytic in $\Sigma_T$.
\end{proof}


\section{Proof of the unique continuation principle} \label{sec3}

In this section we give the proof of the unique continuation principle of solutions to the dual system \eqref{ACP-Dual}.
\begin{proof}[\bf Proof of Theorem \ref{pro-uni-con}]
Let $0\le\nu\le 1$, $1<\mu\le2$,  $\gamma=1/\mu$ or $\gamma=1/2$,  $v_0\in V_\gamma$, $v_1\in L^2(\Omega)$ and let $\omega\subset\Omega$ be an arbitrary non-empty open set. Let $v$ be the unique weak solution to the system \eqref{ACP-Dual} and assume that $v=0$ in $\omega\times (0,T)$.
Since $v=0$ in $\omega\times (0,T)$ and $v:[0,T)\to L^2(\Omega)$ can be analytically extended to the half-plane $\Sigma_T$ (by Theorem \ref{Dual-theo-weak}), it follows that for a.e. $x\in\omega$ and $t\in (-\infty,T)$, 
\begin{align}\label{Exp-v}
v(x,t)=&\sum_{n=1}^\infty (v_0,\varphi_n)_{L^2(\Om)}(T-t)^{\nu(\mu-2)}E_{\mu,\nu(\mu-2)+1}(-\lambda_n(T-t)^\mu)\varphi_n(x)\\&\;\;\;+\sum_{n=1}^\infty (v_1,\varphi_n)_{L^2(\Om)}(T-t)^{\nu(\mu-2)+1}E_{\mu,\nu(\mu-2)+2}(-\lambda_n(T-t)^\mu)\varphi_n(x)=0.
\end{align}
Let $\{\lambda_k\}_{k\in\NN}$ be the set of all eigenvalues of the operator $A$. Let $\{\psi_{k_j}\}_{1\le j\le m_k}\subset L^2(\Omega)$ be the orthonormal basis for $\mbox{Ker}(\lambda_k-A)$. Then \eqref{Exp-v} can be rewritten as
\begin{align}\label{e156}
v(x,t)=&\sum_{k=1}^\infty \left(\sum_{j=1}^{m_k}(v_0,\psi_{k_j})_{L^2(\Om)}\psi_{k_j}(x)\right)(T-t)^{\nu(\mu-2)}E_{\mu,\nu(\mu-2)+1}(-\lambda_k(T-t)^\mu)\notag\\
 &+\sum_{k=1}^\infty \left(\sum_{j=1}^{m_k}(v_1,\psi_{k_j})_{L^2(\Om)}\psi_{k_j}(x)\right)(T-t)^{\nu(\mu-2)+1}E_{\mu,\nu(\mu-2)+2}(-\lambda_k(T-t)^\mu))\notag\\
 &=0,\;\;\;\;\; x\in\omega,\;t\in (-\infty,T).
\end{align}
Let $z\in\CC$ with $\mbox{Re}(z):=\eta>0$ and let $N\in\NN$. Since the system $\{\psi_{k_j}\}\subset L^2(\Omega)$, for $1\le j\le m_k$, $1\le k\le N$ is orthonormal, we have that there is a constant $C>0$ such that
\begin{align*}
&\left\Vert \sum_{k=1}^N \left(\sum_{j=1}^{m_k}(v_0,\psi_{k_j})_{L^2(\Om)}\psi_{k_j}(x)\right)e^{z(t-T)}(T-t)^{\nu(\mu-2)}E_{\mu,\nu(\mu-2)+1}(-\lambda_k(T-t)^\mu))\right\Vert_{L^2(\Omega)}^2\\
\le &\sum_{k=1}^\infty \left(\sum_{j=1}^{m_k}|(v_0,\psi_{k_j})_{L^2(\Om)}|^2\right)e^{2\eta(t-T)}|(T-t)^{\nu(\mu-2)}E_{\mu,\nu(\mu-2)+1}(-\lambda_k(T-t)^\mu)|^2\\
\le & Ce^{2\eta(t-T)}(T-t)^{\nu(\mu-2)}\|v_0\|_{V_\gamma}^2,
\end{align*}
and
\begin{align*}
 &\left\Vert\sum_{k=1}^N\left(\sum_{j=1}^{m_k}(v_1,\psi_{k_j})_{L^2(\Om)}\psi_{k_j}(x)\right)e^{z(t-T)}(T-t)^{\nu(\mu-2)+1}E_{\mu,\nu(\mu-2)+2}(-\lambda_k(T-t)^\mu))\right\Vert_{L^2(\Omega)}^2\\
 \le& Ce^{2\eta(t-T)}(T-t)^{\nu(\mu-2)+1}\|v_1\|_{L^2(\Omega)}^2.
\end{align*}
Define $w_N(\cdot,t)$ by
\begin{align*}
w_N(\cdot,t):=&\sum_{k=1}^N \left(\sum_{j=1}^{m_k}(v_0,\psi_{k_j})_{L^2(\Om)}\psi_{k_j}(x)\right)e^{z(t-T)}(T-t)^{\nu(\mu-2)}E_{\mu,\nu(\mu-2)+1}(-\lambda_k(T-t)^\mu))\\
&+\sum_{k=1}^N\left(\sum_{j=1}^{m_k}(v_1,\psi_{k_j})_{L^2(\Om)}\psi_{k_j}(x)\right)e^{z(t-T)}(T-t)^{\nu(\mu-2)+1}E_{\mu,\nu(\mu-2)+2}(-\lambda_k(T-t)^\mu)).
\end{align*}
We have shown that
\begin{align}\label{nporm-2}
&\|w_N(\cdot,t)\|_{L^2(\Omega)}\notag\\
\le&\left\Vert \sum_{k=1}^N \left(\sum_{j=1}^{m_k}(v_0,\psi_{k_j})_{L^2(\Om)}\psi_{k_j}(x)\right)e^{z(t-T)}(T-t)^{\nu(\mu-2)}E_{\mu,\nu(\mu-2)+1}(-\lambda_k(T-t)^\mu))\right\Vert_{L^2(\Omega)}\notag\\
&+\left\Vert\sum_{k=1}^N \left(\sum_{j=1}^{m_k}(v_1,\psi_{k_j})_{L^2(\Om)}\psi_{k_j}(x)\right)e^{z(t-T)}(T-t)^{\nu(\mu-2)+1}E_{\mu,\nu(\mu-2)+2}(-\lambda_k(T-t)^\mu))\right\Vert_{L^2(\Omega)}\notag\\
&\le Ce^{\eta(t-T)}\left[(T-t)^{\nu(\mu-2)}\|v_0\|_{V_\gamma}+(T-t)^{\nu(\mu-2)+1}\|v_1\|_{L^2(\Omega)}\right].
\end{align}
The right hand side of \eqref{nporm-2} is integrable over $(-\infty,T)$. More precisely, we have that
\begin{align*}
&\int_{-\infty}^Te^{\eta(t-T)}\left[(T-t)^{\nu(\mu-2)}\|v_0\|_{V_\gamma}+(T-t)^{\nu(\mu-2)+1}\|v_1\|_{L^2(\Omega)}\right]\;dt\\
=&\int_0^\infty e^{-\tau}\frac{\tau^{\nu(\mu-2)}}{\eta^{\nu(\mu-2)+1}}\;d\tau |v_0\|_{V_\gamma}+\int_0^\infty e^{-\tau}\frac{\tau^{\nu(\mu-2)+1}}{\eta^{\nu(\mu-2)+2}}\;d\tau |v_1\|_{L^2(\Omega)}\\
=&\frac{\Gamma(\nu(\mu-2)+1)}{\eta^{\nu(\mu-2)+1}}\|v_0\|_{V_\gamma}+\frac{\Gamma(\nu(\mu-2)+2)}{\eta^{\nu(\mu-2)+2}}\|v_1\|_{L^2(\Omega)}.
\end{align*}
By the Lebesgue dominated convergence theorem, we can deduce that for all $\mbox{Re}(z)>0$ and a.e $x\in\Omega$,
\begin{align}\label{e158}
&\int_{-\infty}^Te^{z(t-T)}\left[\sum_{k=1}^\infty \left(\sum_{j=1}^{m_k}(v_0,\psi_{k_j})_{L^2(\Omega)}\psi_{k_j}(x)\right)(T-t)^{\nu(\mu-2)}E_{\mu,\nu(\mu-2)+1}(-\lambda_k(T-t)^\mu)\right.\notag\\
&+\left.\sum_{k=1}^\infty \left(\sum_{j=1}^{m_k}(v_1,\psi_{k_j})_{L^2(\Omega)}\psi_{k_j}(x)\right)(T-t)^{\nu(\mu-2)+1}E_{\mu,\nu(\mu-2)+2}(-\lambda_k(T-t)^\mu)\right]\;dt\notag\\
&=\sum_{k=1}^\infty \sum_{j=1}^{m_k}(v_0,\psi_{k_j})_{L^2(\Omega)}\psi_{k_j}(x)\left(\int_{-\infty}^Te^{z(t-T)}(T-t)^{\nu(\mu-2)}E_{\mu,\nu(\mu-2)+1}(-\lambda_k(T-t)^\mu)dt\right)\notag\\
&+\sum_{k=1}^\infty \sum_{j=1}^{m_k}(v_1,\psi_{k_j})_{L^2(\Omega)}\psi_{k_j}(x)\left(\int_{-\infty}^Te^{z(t-T)}(T-t)^{\nu(\mu-2)+1}E_{\mu,\nu(\mu-2)+2}(-\lambda_k(T-t)^\mu)dt\right)\notag\\
=&\sum_{k=1}^\infty \sum_{j=1}^{m_k}\frac{(v_0,\psi_{k_j})_{L^2(\Omega)}z^{(1-\nu)(\mu-2)+1}}{z^{\mu}+\lambda_k}\psi_{k_j}(x) + \sum_{k=1}^\infty \sum_{j=1}^{m_k}\frac{(v_1,\psi_{k_j})_{L^2(\Omega)}z^{(1-\nu)(\mu-2)}}{z^\mu+\lambda_k}\psi_{k_j}(x),\notag\\
=&\sum_{k=1}^\infty \sum_{j=1}^{m_k}\left(\frac{(v_0,\psi_{k_j})_{L^2(\Omega)}z^{(\nu-1)(\mu-2)+1}}{z^{\mu}+\lambda_k}+\frac{(v_1,\psi_{k_j})_{L^2(\Omega)}z^{(1-\nu)(\mu-2)}}{z^\mu+\lambda_k}\right)\psi_{k_j}(x).
\end{align}
To arrive at \eqref{e158}, we have used the fact that
\begin{align*}
&\int_{-\infty}^Te^{z(t-T)}(T-t)^{\nu(\mu-2)}E_{\mu,\nu(\mu-2)+1}(-\lambda_k(T-t)^\mu)dt\\
=&\int_0^\infty e^{-z\tau}\tau^{\nu(\mu-2)}E_{\mu,\nu(\mu-2)+1}(-\lambda_k\tau^\mu)\;d\tau=\frac{z^{(1-\nu)(\mu-2)+1}}{z^\mu+\lambda_k}
\end{align*}
and
\begin{align*}
&\int_{-\infty}^Te^{z(t-T)}(T-t)^{\nu(\mu-2)+1}E_{\mu,\mu(\mu-2)+2}(-\lambda_k(T-t)^\mu)dt\\
=&\int_0^\infty e^{-z\tau}\tau^{\nu(\mu-2)+1}E_{\mu,\nu(\mu-2)+2}(-\lambda_k\tau^\mu)\;d\tau=\frac{z^{(1-\nu)(\mu-2)}}{z^\alpha+\lambda_k}.
\end{align*}
These identities follow from a simple change of variable and \eqref{lap-ml}.
It follows from \eqref{e156} and \eqref{e158} that
\begin{align*}
\sum_{k=1}^\infty \sum_{j=1}^{m_k}\left(\frac{(v_0,\psi_{k_j})_{L^2(\Omega)}z^{(1-\nu)(\mu-2)+1}}{z^{\mu}+\lambda_k}+\frac{(v_1,\psi_{k_j})_{L^2(\Omega)}z^{(1-\nu)(\mu-2)}}{z^\mu+\lambda_k}\right)\psi_{k_j}(x)=0,\;\;x\in\omega,\;\mbox{Re}(z)>0.
\end{align*}
Letting $\eta:=z^\mu$, we have shown that
\begin{align}\label{e416}
\sum_{k=1}^\infty \sum_{j=1}^{m_k}\left(\frac{(v_0,\psi_{k_j})_{L^2(\Omega)}\eta^{\frac{(1-\nu)(\mu-2)+1}{\mu}}}{\eta+\lambda_k}+\frac{(v_1,\psi_{k_j})_{L^2(\Omega)}\eta^{\frac{(1-\nu)(\mu-2)}{\mu}}}{\eta+\lambda_k}\right)\psi_{k_j}(x)=0,\;\;x\in\omega,\;\mbox{Re}(\eta)>0.
\end{align}
Using  analytic continuation in $\eta$,  we have that the identity  \eqref{e416} holds for every $\eta\in\CC\setminus\{-\lambda_k\}_{k\in\NN}$. Taking a suitable small circle about $-\lambda_l$ and not including $\{-\lambda_k\}_{k\ne l}$ and integrating \eqref{e416} over that circle we get 
\begin{align*}
w_l:= \sum_{j=1}^{m_l}\left[(v_0,\psi_{l_j})_{L^2(\Omega)}(-\lambda_l)^{\frac{(1-\nu)(\mu-2)+1}{\mu}}+(v_1,\psi_{l_j})_{L^2(\Omega)}(-\lambda_l)^{\frac{(1-\nu)(\mu-2)}{\mu}}\right]\psi_{l_j}(x) =0,\;\;x\in\omega.
\end{align*}
Since 
$$(A-\lambda_l)w_l=0\;\mbox{ in }\; \Omega, \;\;w_l=0\;\mbox{ in }\; \omega,$$
and by assumption,  the operator $A$ has the unique continuation property (see \eqref{ucp}), it follows that $w_l=0$ in $\Omega$ for every $l$. Since $\{\psi_{l_j}\}_{1\le j\le m_k}$ is linearly independent in $L^2(\Omega)$, we obtain that 
\begin{align*}
(v_0,\psi_{l_j})_{L^2(\Om)}(-\lambda_l)^{\frac{\beta+1}{\mu}}+(v_1,\psi_{l_j})_{L^2(\Om)}(-\lambda_l)^{\frac{\beta}{\mu}}=0\;\;\mbox{for}\;1\le j\le m_k,\;k\in\NN\;\mbox{and}\;\beta=(1-\nu)(\mu-2).
\end{align*}
This implies that
\begin{align}\label{zero}
0=&(-\lambda_l)^{\frac{\beta+1}{\mu}}(v_0,\psi_{l_j})_{L^2(\Om)}+(-\lambda_l)^{\frac{\beta}{\mu}}(v_1,\psi_{l_j})_{L^2(\Om)}\notag\\
=&\lambda_l^{\frac{\beta+1}{\mu}}\left[\cos(\frac{\beta+1}{\mu}\pi)+i\sin(\frac{\beta+1}{\mu}\pi)\right](v_0,\psi_{l_j})_{L^2(\Om)}+\lambda_l^{\frac{\beta}{\mu}}\left[\cos(\frac{\beta}{\mu}\pi)+i\sin(\frac{\beta}{\mu}\pi)\right](v_1,\psi_{l_j})_{L^2(\Om)}.
\end{align}
It follows from \eqref{zero} that $v_0=0=v_1$ in $\Omega$.  Hence, $v=0$ in $\Omega\times (0,T)$ and the proof is finished.
\end{proof}

\section{Proof of the controllability result}\label{prof-ma-re}\label{sec4}

Before proceeding with the proof, we show that to study the memory approximate controllability of  \eqref{main-EQ}, it suffices to consider the case $u_0=u_1=0$.

\begin{remark} \label{rem-41}
{\em
Consider the following two systems:
\begin{align}\label{EQ-0}
\begin{cases}
\mathbb D_t^{\mu,\nu} v+Av=f|_{\omega}\;\;&\mbox{ in }\; \Omega\times (0,T),\\
\left(I_t^{(1-\nu)(2-\mu)}v\right)(\cdot,0)=0 &\mbox{ in }\;\Omega,\\
\left(\partial_tI_t^{(1-\nu)(2-\mu)}v\right)(\cdot,0)=0 &\mbox{ in }\;\Omega,
\end{cases}
\end{align}
and
\begin{align}\label{EQ-U0}
\begin{cases}
\mathbb D_t^{\mu,\nu} w+Aw=0\;\;&\mbox{ in }\; \Omega\times (0,T),\\
\left(I_t^{(1-\nu)(2-\mu)}w\right)(\cdot,0)=u_0 &\mbox{ in }\;\Omega,\\
\left(\partial_tI_t^{(1-\nu)(2-\mu)}w\right)(\cdot,0)=u_1 &\mbox{ in }\;\Omega.
\end{cases}
\end{align}
Given $u_0\in V_\gamma$ and $u_1\in L^2(\Om)$, let  $w$ be the weak solution of  \eqref{EQ-U0}. Assume that  \eqref{EQ-0} is memory  approximately controllable. Then, for every $\varepsilon>0$ and $(\tilde u,\tilde v)\in V_\gamma\times L^2(\Om)$, there exists a control function $f\in L^p((0,T);L^2(\omega))$ such that the corresponding unique weak solution $v$ of  \eqref{EQ-0} satisfies
\begin{align}\label{control}
&\left\|I_t^{(1-\nu)(2-\mu)}v(\cdot,T)-(\tilde u-I_t^{(1-\nu)(2-\mu)}w(\cdot,T))\right\|_{V_\gamma}\notag\\
+&\left\|\partial_tI_t^{((1-\nu)(2-\mu)}v(\cdot,T)-(\tilde v-\partial_tI_t^{((1-\nu)(2-\mu)}w(\cdot,T))\right\|_{L^2(\Omega)}\le\varepsilon.
\end{align}
By definition, we have that the function $v+w$ solves the system \eqref{main-EQ}, and it follows from \eqref{control} that  
\begin{align*}
\|I_t^{(1-\nu)(2-\mu)}(v+w)(\cdot,T)-\tilde u\|_{V_\gamma}+\|\partial_tI_t^{(1-\nu)(2-\mu)}(v+w)(\cdot,T)-\tilde v\|_{L^2(\Omega)}\le\varepsilon.
\end{align*}
Hence, \eqref{main-EQ} is memory approximately controllable.   Therefore,  in our study we consider the system \eqref{main-EQ} with $u_0=u_1=0$.
}
\end{remark}

\begin{proof}[\bf Proof of Theorem \ref{second-Theo}]
Let $0\le\nu\le 1$, $1<\mu\le 2$, $\gamma=1/\mu$ and $p>1/(\mu-1)$ (or $\gamma=1/2$ and $p>2/\mu$) and $f\in L^p((0,T);L^2(\omega))$.
Assume that the operator $A$ has the unique continuation property in the sense of \eqref{ucp}.
Let $u$ be the unique weak solution of \eqref{main-EQ} with $u_0=u_1=0$ and $v$ the unique weak solution of the adjoint system \eqref{ACP-Dual} with $v_0\in V_\gamma$ and $v_1\in L^2(\Omega)$. 

Firstly, it follows from \eqref{norm-V} that $v\in L^q((0,T);V_\gamma)\hookrightarrow L^q((0,T);L^2(\Omega))$ for every $1\le q< 1/\nu(2-\mu)$.  Secondly, let $p^\star$ denote the conjugate exponent of $p$. Since by assumption $p>1/(\mu-1)$, we have that $1\le p^\star=1+1/(p-1) < 1/(2-\mu)$. Since $\nu(2-\mu)\le 2-\mu$, we have that $v\in L^{p^\star}((0,T);L^2(\Omega)$.
In addition, it follows from the proofs of Theorem \ref{theo-weak} and Theorem \ref{Dual-theo-weak} that $u$ and $v$ have the required regularity to apply the integration by parts formula \eqref{IP01}.
Then, integrating by parts by using \eqref{IP01}  together with the fact that for every $u,v\in V_{\gamma}$,
 $$\langle Au,v\rangle_{V_{-\gamma},V_\gamma}=\langle u,Av\rangle_{V_{\gamma},V_{-\gamma}}$$ 
we get
\begin{align*}
0=&\int_0^{T}\langle \mathbb D_t^{\mu,\nu} u+Au-f,v\rangle_{V_{-\gamma},V_\gamma}\;dt\notag\\
=&\int_0^{T}\langle v,\mathbb D_t^{\mu,\nu} u+Au\rangle_{V_{\gamma},V_{-\gamma}}\; dt -\int_0^{T}\int_{\omega}fv\;dxdt\notag\\
=&\int_0^{T}\langle\mathbb D_{t,T}^{\mu,(1-\nu)} v+Av,u\rangle_{V_{-\gamma},V_\gamma}\;dt -\int_0^{T}\int_{\omega}fv\;dxdt\notag\\
&+\int_{\Omega}\left[\partial_tI_t^{(1-\nu)(2-\mu)}u(x,T) I_{t,T}^{\nu(2-\mu)}v(x,T)+I_t^{(1-\nu)(2-\mu)}u(x,T)D_{t,T}^{1-\nu(2-\mu)}v(x,T)\right]\;dx\notag\\
=&\int_{\Omega}\left[\partial_tI_t^{(1-\nu)(2-\mu)}u(x,T) I_{t,T}^{\nu(2-\mu)}v(x,T)+I_t^{(1-\nu)(2-\mu)}u(x,T)D_{t,T}^{1-\nu(2-\mu)}v(x,T)\right]\;dx\notag\\
&-\int_0^{T}\int_{\omega}fv\;dxdt.
\end{align*}
We  have shown  that
\begin{align}\label{eq41}
\int_{\Omega}\left[\partial_tI_t^{(1-\nu)(2-\mu)}u(x,T)v_0+I_t^{(1-\nu)(2-\mu)}u(x,T)v_1\right]\;dx=\int_0^T\int_{\omega}fv\;dxdt.
\end{align}
To prove that the set 
$$\left\{\left( I_t^{(1-\nu)(2-\mu)}u(\cdot,T),\partial_t I_t^{(1-\nu)(2-\mu)}u(\cdot,T)\right):\; f\in L^p((0,T);L^2(\omega))\right\}$$ 
is dense in $V_\gamma\times L^2(\Omega)$, we have to show that if $(v_0,v_1)\in V_\gamma\times L^2(\Omega)$ is such that
\begin{align}\label{eq42}
\int_{\Omega}\left[\partial_tI_t^{(1-\nu)(2-\mu)}u(x,T)v_0(x)+I_t^{(1-\nu)(2-\mu)}u(x,T)v_1(x)\right]\;dx=0
\end{align}
for every $f\in L^p((0,T);L^2(\omega))$, then $v_0=v_1=0$. Indeed, let $v_0$ and $v_1$ satisfy \eqref{eq42}. It follows from \eqref{eq41} and \eqref{eq42} that
\begin{align*}
\int_0^T\int_{\omega}fv\;dxdt=0
\end{align*}
for every $f\in L^p((0,T);L^2(\omega))$.  By the fundamental lemma of the calculus of variations, we have that
\begin{align*}
v=0 \;\mbox{ in }\; \omega\times (0,T).
\end{align*}
It follows from Theorem \ref{pro-uni-con} that
\begin{align*}
v=0  \;\mbox{ in }\; \Omega\times (0,T).
\end{align*}
Since the solution $v$ of \eqref{ACP-Dual} is unique,  we have that $v_0=v_1=0$ on $\Omega$. The proof  is finished.
\end{proof}

We conclude the paper with the following observation.

\begin{remark}\label{rem-map-ucp}
{\em 
$0\le\nu\le 1$, $1<\mu\le 2$, $\gamma=1/\mu$ and $p>1/(\mu-1)$ (or $\gamma=1/2$ and $p>2/\mu$),  $f\in L^p((0,T);L^2(\omega))$ and consider the following mapping:
\begin{align*}
F: L^p((0,T), L^2(\omega))\to V_\gamma\times L^2(\Omega),\; f\mapsto\Big( I_t^{(1-\nu)(2-\mu}u(\cdot,T),\partial_t I_t^{(1-\nu)(2-\mu}u(\cdot,T)\Big),
\end{align*}
where $u$ is the unique weak solution of \eqref{main-EQ} with $u_0=u_1=0$.  It is easy to see that the system \eqref{main-EQ} is memory approximately controllable in time $T>0$ if and only if the range of $F$, that is, $\mbox{Ran}(F)$ is dense in $V_\gamma\times L^2(\Omega)$. This is equivalent to $\mbox{Ker}(F^\star)=\{(0,0)\}$, where $F^\star$ is the adjoint of $F$. It follows from the proof of Theorem \ref{second-Theo} (more precisely from \eqref{eq41}) that $F^\star$ is the mapping given by
\begin{align*}
F^\star: V_\gamma\times L^2(\Omega)\to L^p((0,T);L^2(\omega)),\;(v_0,v_1)\mapsto v\big|_{\omega\times(0,T)},
\end{align*}
where $v$ is the unique solution of the adjoint system \eqref{ACP-Dual}. Again $\mbox{Ker}(F^\star)=\{(0,0)\}$ is equivalent to the unique continuation principle, namely,
\begin{align*}
\Big(v\;\mbox{ solution of }\; \eqref{ACP-Dual}:\;  v\big|_{\omega\times(0,T)}=0\Big) \Longrightarrow \;v_0=v_1=0\;\mbox{ in }\;\Omega \Longrightarrow v=0 \mbox{ in } \Omega\times (0,T).
\end{align*}
}
\end{remark}

 \bibliographystyle{plain}
\bibliography{biblio}

\end{document}